\documentclass[journal]{IEEEtran}
\usepackage{cite}

  \usepackage[pdftex]{graphicx}
\usepackage[cmex10]{amsmath}
\usepackage{array}
\usepackage{fixltx2e}
\usepackage{dblfloatfix}

\ifCLASSOPTIONcaptionsoff
 \usepackage[nomarkers]{endfloat}
 \let\MYoriglatexcaption\caption
 \renewcommand{\caption}[2][\relax]{\MYoriglatexcaption[#2]{#2}}
\fi
\usepackage{amsfonts}
\usepackage{amssymb}
\usepackage{accents}

\newcommand\Pmi{P_{\mathrm{m}i}}

\newcommand\Pmo{P_{\mathrm{m}1}}
\newcommand{\argmax}{\operatornamewithlimits{arg\,max}}
\newcommand\Pmt{P_{\mathrm{m}2}}
\newcommand\delhh{\mathit{\Delta\mathcal{H}}_\mathrm{max}}
\newcommand\dspre{\boldsymbol{\delta}^s_{\mathrm{pre}}}
\usepackage{accents}
\newlength{\dhatheight}

\begin{document}
%
\title{An analytical critical clearing time for parametric analysis of transient stability in power systems}
%
%
%

\author{Lewis~Roberts,~\IEEEmembership{Student~Member,~IEEE,}
		Alan~Champneys,
        Keith~Bell,~\IEEEmembership{Member,~IEEE,}
        Mario~di~Bernardo,~\IEEEmembership{Fellow,~IEEE.}
\thanks{L. Roberts, A. Champneys and M. di Bernardo are in the Department of Engineering Mathematics, University of Bristol, UK. e-mail: lewis.roberts@bristol.ac.uk. M. di Bernardo is also with the Department of Electrical Engineering and ICT, University of Naples Federico II, Italy.}
\thanks{K. Bell is the ScottishPower Professor of Smart Grids at the University of Strathclyde, Glasgow UK.}
}

%
%

\markboth{}%
{Shell \MakeLowercase{\textit{et al.}}: Bare Demo of IEEEtran.cls for Journals}
%



\maketitle

\begin{abstract}
An analytic approximation for the critical clearing time (CCT) metric is derived from direct methods for power system stability. The formula has been designed to incorporate as many features of transient stability analysis as possible such as different fault locations and different post-fault network states. The purpose of this metric is to analyse trends in stability (in terms of CCT) of power systems under the variation of a system parameter. We demonstrate the performance of this metric to measure stability trends on an aggregated power network, the so-called two machine infinite bus network, by varying load parameters in the full bus admittance matrix using numerical continuation. Our metric is compared to two other expressions for the CCT which incorporate additional non-linearities present in the model.
\end{abstract}

%
\IEEEpeerreviewmaketitle

\section{Introduction}
%
%
%
%

\IEEEPARstart{T}{he} complex dynamics of electric power systems have long been the subject of intense research particularly in the area of stability. Effective stability metrics provide control inputs and assist the system operator to ensure that a power system maintains synchrony after the network suffers a fault, i.e. that it exhibits transient stability. A traditional transient stability metric for short circuit faults on a power network is the so-called \textit{critical clearing time} (CCT) \cite{Kundur1994,Machowski2008book}. The CCT provides an upper bound on the duration of a short circuit on a power network before it is removed - `cleared' - by the action of protection mechanisms to isolate the faulted circuit such that the system will regain synchronisation once the fault is cleared. In general, the CCT is a useful metric for power system design; by allowing the severity of different situations and the effectiveness of different interventions (generation re-dispatches, control modifications or network reinforcements) to be compared.

Currently, there are practical developments in power systems that promise to radically change power system dynamic behaviour. For example, the gradual substitution of power generated from large, synchronous machines by asynchronous machines or power fed via power electronic interfaces (e.g. wind farms, solar PV and HVDC interconnections to other systems), in addition to the changing nature of electrical loads \cite{yamashita2012modelling}. Previous work in the literature \cite{Slootweg2002} has investigated the effect of changing loads on system stability by repeating fault studies for different loading levels. As a consequence, there is value in articulating metrics that exploit theoretical, if simplified, descriptions of the system which can provide a deep understanding of the impact of a wide range of features of the network from parametric investigations. This can inform efforts to design strategies to mitigate possible instabilities in the system.  

In the recent literature, alternative methodologies have been used to study stability when modelling a power system using the so-called swing equations \cite{Kundur1994,Pai1989book}. These include synchronisation \cite{BULLO2012}, non-linear dynamics \cite{Susuki2011}, bifurcation theory \cite{ajjarapu1992bifurcation}, passivity-based methods \cite{Galaz2003} and the computation of basins of attraction \cite{Hasegawa1999}. Direct methods \cite{Chiang2011book} cast the swing equations in an energetic framework to provide a critical energy boundary for the whole system during a fault. Despite the difficulty of including non-negligible transfer conductances in direct methods \cite{Sastry1973}, their advantages include the possible estimation of an analytical stability boundary and relatively quick computation. Also, they require no need for further simplifications of a power system beyond the swing equation model and they can be applied to any system that can be parametrised. The system operator can use this analytical stability metric for initial safety checks and to assess the stability margins of the system once a fault has been cleared. 

One of the drawbacks of the direct methods is that it is difficult to predict when the system energy will cross the critical energy boundary because of the non-linear nature of the system dynamics. So-called fault trajectory sensitivity techniques \cite{Laufenberg1997,Fouad1981,Hiskens2000} have been proposed to consider the effect of parameters on stability by linearising about the trajectory of a fault in state space with respect to a given parameter. Furthermore, a method for computing a so-called ``direct CCT'' has been proposed \cite{Nguyen2003} which is based on linearising the power system model about a specific fault trajectory with respect to the system energy itself. An estimate of the CCT is then found by extrapolation. However, to our knowledge, an analytic CCT metric is only available for induction generators \cite{Grilo2007} and there is no analytic estimate of the CCT developed for a network of synchronous generators. 

The aim of this paper is to propose a new analytic expression of the CCT. This estimate is derived by recasting the energetic metric used in the direct methods in terms of a metric in time by simplifying the energy functions and the dynamics during a fault. As is true for direct methods in general, our metric can serve as a lower bound to the true CCT for lossless power systems (or for power networks with small transfer conductances \cite{Chiang1995}) and can be applied to systems suffering a large fault at any location on a network. However, the purpose of our metric is to capture trends in stability as network parameters are varied and as such, the investigations in this paper are limited to aggregated or clustered power networks. (See \cite[Chapter 14]{Machowski2008book} for aggregation techniques.)

In general, a power network's topology changes from its original structure when a fault is cleared. This is generally due to some switching action that isolates the region of the network that suffers the fault. Choosing the best strategy to quickly identify a need for and carry out this action is a crucial step in maintaining the stability and synchronisation of a power system. There is some uncertainty regarding the success and speed of protection actions, and, as a consequence, power flows may need to be restricted and more expensive, or higher carbon, power sources utilised. We argue that choices both in operation and the design of the system and its control can be facilitated by parameter investigations of power system stability models such as the swing equation and applying quick but effective stability metrics to illustrate the effect of a given parameter value change. A rigorous study of the strategies available to the system operator could be provided in-part by the continuous variation of model parameters, which could possibly uncover optimal parameter values to maximise stability at the design stage or on-line. The analytic CCT metric derived in this paper is able to capture sensitivities in stability of a given fault as a network parameter is varied. In particular, this paper studies the effect of a load parameter on the stability of a given fault in an aggregated network.

The rest of this paper is organised as follows: In Section~\ref{sec:fault} we formulate a CCT estimate denoted $\tau_\mathrm{H}$ (where the subscript `$\mathrm{H}$' signifies `Hamiltonian') using direct methods and introduce the aggregate network used to conduct our investigations, the two-machine infinite bus (TMIB) network. By considering polynomial approximations of $\tau_\mathrm{H}$ we derive an analytic CCT metric denoted $\tau_\mathrm{A}$ (where the subscript `$\mathrm{A}$' signifies `analytic') in Section~\ref{sec:low_bound}. A parametric investigation of the effect on stability of different loadings on an aggregated network given a particular fault is presented in Section~\ref{sec:3genExample} and finally, conclusions are drawn in Section~\ref{sec:discussion} together with suggestions for future work.

\section{Fault analysis using energy functions}
\label{sec:fault}

\subsection{Model description}
We consider the classic swing equation model \cite{Machowski2008book,Kundur1994} to describe the stability effects of transient faults on a power system with synchronous generation. The generators are modelled as voltage sources behind reactances and the loads on the network are of constant impedance. In general, generators have small losses due to damping \cite{Anderson2002book} so without loss of generality we assume zero damping for generators. This model can be written as a set of coupled one-dimensional ordinary differential equations (ODEs), which describe the dynamics of the rotor angles of each synchronous generator $i \in \{1,\dots,n\}$  in a network by considering Newton's second law of dynamics. In vector form the equation is
\begin{equation}
\dot{\mathbf{x}} = \mathbf{F}(\mathbf{x}),
\label{eq:genMMSE}
\end{equation}
where
\begin{equation*}
\mathbf{x} = 
\begin{pmatrix}
\boldsymbol{\delta} \\
\boldsymbol{\omega}
\end{pmatrix}
\end{equation*}
and
\begin{equation*}
\mathbf{F}\left(\begin{pmatrix}
\boldsymbol{\delta} \\
\boldsymbol{\omega}
\end{pmatrix}\right) = 
\begin{pmatrix}
\boldsymbol{\omega} \\
\mathbf{A}(\boldsymbol{\delta})
\end{pmatrix}.
\end{equation*}
The vectors $\boldsymbol{\delta} = [\delta_1,\dots,\delta_n]^T$ and $\boldsymbol{\omega} = [\omega_1,\dots,\omega_n]^T$ are the generator rotor angles and angular speeds respectively, and the elements of the vector function $\mathbf{A}(\boldsymbol{\delta})$ are
\begin{equation*}
A_i(\boldsymbol{\delta}) = \frac{1}{M_i}\left( \Pmi - P_{\mathrm{e}i}(\boldsymbol{\delta}) \right),
\end{equation*}
where $M_i = \frac{2H_i}{\omega_0}$ is a lumped parameter,  $\omega_0 = 2\pi f$ (where $f$ is the grid frequency: $50\,\mathrm{Hz}$ in Europe), $H_i$ is the inertia constant, $\Pmi$ is the mechanical power input and $P_{\mathrm{e}i}(\boldsymbol{\delta})$ is the electrical power output.

The loads on the power system are assumed to be constant impedance loads such that Kron reduction \cite{Dorfler2013a} can be applied to the network. Therefore, the swing equations describe the dynamics of a reduced network comprising of constant voltage sources connected through a network of impedances \cite{Machowski2008book}. The total power consumed by conductive loads at generator $i$ is given by
\begin{equation}
P_i(\boldsymbol{\delta}) = E_i^2 G_{ii} +\sum_{k\neq i}^n|E_i||E_k|G_{ik}\cos(\delta_i - \delta_k),
\label{eq:cond_elec}
\end{equation} 
where $E_i = |E_i|e^{j\delta_i}$ is the internal voltage of generator $i$ ($|E_i|$ assumed constant), $G_{ik}$ is the conductance between generators $i$ and $k$ and $G_{ii}$ is the shunt conductance at bus $i$. The total electric power leaving generator $i$ is
\begin{equation}
P_{\mathrm{e}i}(\boldsymbol{\delta}) = P_i(\boldsymbol{\delta}) + \sum_{k\neq i}^n\bar{P}_{ik}\sin(\delta_i - \delta_k),
\label{eq:elec_power_term}
\end{equation} 
where $\bar{P}_{ik}=|E_i||E_i|B_{ik}$ is the maximum active power flow between generators $i$ and $k$ and $B_{ik}$ is the susceptance of the network connection between node $i$ and node $k$.  
The admittances $Y_{ik} = G_{ik} + jB_{ik}$ are the elements of the (symmetric) reduced bus admittance matrix $\mathbf{Y}_{\mathrm{red}} \in \mathbb{C}^{n \times n}$. Kron reduction is fundamentally a matrix operation permitted by applying Kirchoff's laws to a power network and constructing $\mathbf{Y}_{\mathrm{red}}$ from a larger bus admittance matrix $\mathbf{Y}_{\mathrm{BUS}} \in \mathbb{C}^{N\times N}$ where  $N \geq 2n$. The bus matrix $\mathbf{Y}_{\mathrm{BUS}}$ is a block matrix which contains the full topology and load distribution (including the synchronous reactance) of a power network with $n$ synchronous generators. 

A stationary point of the system (\ref{eq:genMMSE}) solves the equation $\mathbf{F}(\mathbf{x}) = \mathbf{0}$ (where $\mathbf{0}$ is a column vector of all zeros) and is denoted $\mathbf{x}^* = [{\boldsymbol{\delta}^*}^T;{\boldsymbol{\omega^*}}^T]^T$. A solution for (\ref{eq:genMMSE}) starting from initial conditions $\mathbf{x}(0)$ is written generically as
\begin{equation}
\mathbf{x}(t) = \boldsymbol{\Phi}(t;\mathbf{x}(0)) , \quad t \geq 0.
\label{eq:sol_swing}
\end{equation}

\subsection{Fault analysis}
\label{sec:fault_anal}
\subsubsection{Stability analysis of transient faults}
The objective of transient fault analysis is to investigate whether a system will remain stable once a fault has been cleared and, ideally, no further action from the system operator would be required. We assume, without loss of generality, that the moment a power system suffers a short-circuit is at time $t=0$ and the fault is cleared at time $t_\mathrm{cl}$. These two points in time define three distinct regimes in order to analyse the dynamics of a fault on a power system. These are $(i) \;t<0$ \textit{(pre-fault)}, $(ii)\; 0 \leq t < t_\mathrm{cl}$ (\textit{fault-on}) and $(iii) \; t \geq t_\mathrm{cl} $ (\textit{post-fault}). 

The fault analysis method in \cite[Chapter 2]{Anderson2002book} (recently summarised in \cite{Gajduk2014a}) for power networks with constant impedance loads, is employed in this paper. Each regime has a different bus matrix $\mathbf{Y}_\mathrm{BUS}$ (and therefore reduced admittance matrix $\mathbf{Y}_{\mathrm{red}}$) which will change the values of the parameters $G_{ii}$, $G_{ik}$ and $B_{ik}\;\mathrm{for\;all}\;i,k$ in the vector function (\ref{eq:genMMSE}). Therefore, three separate sets of equations of the form (\ref{eq:genMMSE}) are required to model the power system for all time given by
\begin{equation}
\dot{\mathbf{x}} = 
\begin{cases}
\mathbf{F}_\mathrm{pre}(\mathbf{x}) \quad t<0 \\
\mathbf{F}_\mathrm{on}(\mathbf{x}) \quad 0 \leq t < t_\mathrm{cl} \\
\mathbf{F}_\mathrm{post}(\mathbf{x}) \quad t \geq t_\mathrm{cl}
\end{cases},
\label{eq:fault_model}
\end{equation}
where the labels `$\mathrm{pre}$', `$\mathrm{on}$' and `$\mathrm{post}$' refer to the parameter values for the system in regimes $(i)$, $(ii)$ and $(iii)$ respectively. Pre-fault, a power system is assumed to be balanced and therefore we assume that (\ref{eq:fault_model}) is located at a stable (`s') equilibrium point 
\begin{equation*}
\mathbf{x}^s_\mathrm{pre} =
\begin{pmatrix}
\dspre \\
\mathbf{0}
\end{pmatrix}
\end{equation*}
for $t<0$ where ${|\delta^s_{\mathrm{pre},i}-\delta^s_{\mathrm{pre},k}| < \pi/2}\;\mathrm{for\;all}\;i,k$. The dynamics for $t\geq 0$ are given by
\begin{equation}
\mathbf{x}_\mathrm{on}(t) = \boldsymbol{\Phi}_\mathrm{on}(t;\mathbf{x}_\mathrm{on}(0) = \mathbf{x}_\mathrm{pre}^s), \quad 0 \leq t < t_\mathrm{cl}
\label{eq:fault_traj}
\end{equation}
during the fault and
\begin{equation}
\mathbf{x}_\mathrm{post}(t) = \boldsymbol{\Phi}_\mathrm{post}(t,\mathbf{x}_\mathrm{post}(0) = \mathbf{x}_\mathrm{on}(t_\mathrm{cl})), \quad t \geq t_\mathrm{cl}
\label{eq:pft}
\end{equation}
after the fault. From these expressions we can define the CCT, denoted $\tau$, formally as the maximum value of $t_\mathrm{cl}$ such that in the post-fault trajectory (\ref{eq:pft}) there is one full swing of the rotor angles before some pairs of rotors angles begin to diverge \cite{Anderson2002book}; this is also known as \textit{first swing} stability and is generally found algorithmically using power systems software packages.

\subsubsection{A CCT approximation using energetic methods}
\label{sec:faults}

In general, a conservative metric for the local stability of systems of the form (\ref{eq:genMMSE}) can be found by constructing a suitable Lyapunov function. Direct methods use so-called energy functions \cite{Pai1989book}, which can also serve as Lyapunov functions, to measure the global stability of such systems. A stability boundary is constructed in terms of a critical system energy $\mathcal{E}_\mathrm{c}$ in the post-fault regime and a power system is classified as unstable when the total system energy surpasses this critical energy. 

The total system energy can be measured when a power system is modelled as a Hamiltonian system. However, the power consumed by the loads $P_i(\boldsymbol{\delta})$ is a path-dependent quantity \cite{Pai1989book} and cannot be modelled exactly by a conservative system. The survey paper \cite{Ribbens-Pavella1985} collects numerous attempts that have been used to approximate this term so that an appropriate Hamiltonian system can be used. The most accepted technique \cite[p. 231]{Machowski2008book} models the power consumed by the loads as a constant term given by 
\begin{equation}
P_{\mathrm{a}i} = P_{i}(\boldsymbol{\delta}^s),
\label{eq:diss_const}
\end{equation}
where the point $\mathbf{x}^s = [{\boldsymbol{\delta}^s}^T,{\mathbf{0}}^T]^T$ is a stable stationary point in the post-fault regime which solves $\mathbf{F}_\mathrm{post}
(\mathbf{x}^s) = \mathbf{0}$ with $|\delta^s_i-\delta^s_k| < \pi/2\;\mathrm{for\;all}\;i,k$. The dynamics of a power system with assumption (\ref{eq:diss_const}) employed can be written as
\begin{equation}
\dot{\mathbf{x}} = \hat{\mathbf{F}}(\mathbf{x},\mathbf{x}^s),
\label{eq:genMMSEham}
\end{equation}
where terms depending on the conductive parts of loads are isolated to obtain the vector function $\hat{\mathbf{F}}(\mathbf{x},\mathbf{x}^s)$. This function has a similar structure as in (\ref{eq:genMMSE}) where
\begin{equation*}
\hat{\mathbf{F}}\left(
\begin{pmatrix}
\boldsymbol{\delta} \\
\boldsymbol{\omega}
\end{pmatrix},
\begin{pmatrix}
\boldsymbol{\delta}^s \\
\mathbf{0}
\end{pmatrix}
\right) = 
\begin{pmatrix}
\boldsymbol{\omega} \\
\hat{\mathbf{A}}(\boldsymbol{\delta},\boldsymbol{\delta}^s)
\end{pmatrix},
\end{equation*}
and the elements of the vector $\hat{\mathbf{A}}(\boldsymbol{\delta},\boldsymbol{\delta}^s)$ are given by
\begin{equation*}
\hat{A}_i(\boldsymbol{\delta},\boldsymbol{\delta}^s) = \frac{1}{M_i}\left( \Pmi - \hat{P}_{\mathrm{e}i}(\boldsymbol{\delta},\boldsymbol{\delta}^s) \right),
\end{equation*}
with
\begin{equation*}
\hat{P}_{\mathrm{e}i}(\boldsymbol{\delta},\boldsymbol{\delta}^s) = P_{i}(\boldsymbol{\delta}^s) + \sum_{k\neq i}^n\bar{P}_{ik}\sin(\delta_i - \delta_k).
\end{equation*}

The Hamiltonian function
\begin{equation}
\mathcal{H}(\mathbf{x}) = \mathcal{E}_\mathrm{kin}(\boldsymbol{\omega}) + \mathcal{E}_\mathrm{pot}(\boldsymbol{\delta}),
\label{eq:ham}
\end{equation}
quantifies the post-fault system energy for a system of the form (\ref{eq:genMMSEham}) and is the sum of the kinetic energy $\mathcal{E}_\mathrm{kin}(\boldsymbol{\omega})$ and the potential energy $\mathcal{E}_\mathrm{pot}(\boldsymbol{\delta})$ for a power system with $n$ generators where
\begin{equation}
\mathcal{E}_\textrm{kin}(\boldsymbol{\omega}) = \sum_{i=1}^n{\frac{1}{2}M_i \omega_i^2},
\label{eq:kin_multi}
\end{equation}
and
\begin{equation}
\mathcal{E}_\textrm{pot}(\boldsymbol{\delta}) = - \sum_{i=1}^n { (\Pmi-P_{\mathrm{a}i}) \delta_i} - \sum_{\substack{i=1 \\ k > i}}^n {\bar{P}_{ik} \cos (\delta_i - \delta_k)}.
\label{eq:pot_multi}
\end{equation}

An approximation of the CCT, denoted $\tau_\mathrm{H}$ can be found by integrating the dynamics of the system during a fault until the system energy reaches the critical boundary $\mathcal{E}_\mathrm{c}$ (which will be computed later). More specifically, such an estimate can be obtained by observing the first instance that the Hamiltonian 
\begin{equation}
\mathcal{H}(\mathbf{x}_\mathrm{on}(t)) = \mathcal{E}_\mathrm{c},
\label{eq:solve_CCT}
\end{equation}
for $t>0$ where, in general, the energy difference
\begin{equation}
\delhh :=  \mathcal{E}_\mathrm{c} - \mathcal{H}(\mathbf{x}^s_\mathrm{pre}),
\end{equation}
is positive for a suitably chosen post-fault network.  Note that, the power system during the fault is not modelled as a Hamiltonian. The CCT approximation $\tau_\mathrm{H}$ is much faster to compute than the traditional CCT because the dynamics of the post-fault system (\ref{eq:pft}) do not need to be computed. 

The so-called \textit{closest UEP} (unstable equilibrium point) method \cite{Chiang2011book} is used to find the critical system energy $\mathcal{E}_\mathrm{c}$ in this work because, although it is the most conservative method (compared to the controlling UEP method or the potential energy boundary surface method \cite{Chiang2011book}) it can be applied to any power system without considering the specific fault that a system suffers. In the presence of large linear loads in the network, the use of direct methods might lead to overestimates of the actual stability boundary \cite{Athay1979} however, the intention is to study the effect of stability trends, so the closest UEP serves as an adequate method to capture the system energy for initial parametric studies.

The critical energy boundary computed by the closest UEP method is defined as
\begin{equation}
\mathcal{E}_\mathrm{c} =  \mathcal{E}_\mathrm{pot}(\boldsymbol{\delta}^u_{c}) = \min\{ \mathcal{E}_\mathrm{pot}(\boldsymbol{\delta}^u_{1}),\dots, \mathcal{E}_\mathrm{pot}(\boldsymbol{\delta}^u_{m})\}.
\label{eq:cuep_def}
\end{equation}
The point $\mathbf{x}^u_{c} \in \mathcal{S}$ is the so-called closest UEP where 
\begin{equation}
\mathcal{S} = \{ \mathbf{x}_1^u,\dots, \mathbf{x}_m^u\}
\label{eq:unst_set}
\end{equation}
is the set of all `type-1' \cite{Chiang2011book} unstable equilibria of (\ref{eq:genMMSEham}) where 
\begin{equation*}
\mathbf{x}_i^u = 
\begin{pmatrix}
\boldsymbol{\delta}_i^u \\
\mathbf{0}
\end{pmatrix}.
\end{equation*}

\subsection{An aggregate network}

In order to study trends in stability under parametric variations, the dynamics of each generator in a network can be grouped into synchronous regions according to the electrical distance between individual generators. Previous studies \cite{Hughes2006,Ulbig,Johnstone2014} have used an aggregate power network model to study the global dynamics of a power system. Typically, the models presented in these references study the GB power network using a three bus network with the inertia of one machine at least two order of magnitudes larger than the other two. These models lend themselves well to be studied using a so-called two machine infinite bus (TMIB) system. This network structure has been previously studied in \cite{Llamas1995,Chiang2011book,Anghel2013,Guo2000}. The ODE in the form (\ref{eq:genMMSEham}) for this system is given by 
\begin{equation}
\begin{aligned}
\dot{\delta}_1 &= \omega_1 \\
\dot{\delta}_2 &= \omega_2 \\
\dot{\omega}_1 &= \frac{1}{M_1}\left[(\Pmo-P_{\mathrm{a}1}) - \bar{P}_{13}\sin(\delta_1) - \bar{P}_{12}\sin(\delta_1-\delta_2)\right] \\
\dot{\omega}_2 &= \frac{1}{M_2}\left[(\Pmt-P_{\mathrm{a}2}) - \bar{P}_{23}\sin(\delta_2) - \bar{P}_{12}\sin(\delta_2-\delta_1)\right] \\
\end{aligned}
\label{eq:TMIBode}
\end{equation}
where we have employed assumption (\ref{eq:diss_const}) to get a conservative system.
After Kron reduction there are three interconnected buses in the network: two buses connected to synchronous generators and an infinite bus (bus 3). The infinite bus models the dynamics of a large section of a network as a generator with infinite inertia and constant internal voltage  $E_3$. As such, $\delta_3$ is a constant and without loss of generality we can set $\delta_3 = 0$ and use it as a reference point for the other two rotor angles. 

\begin{figure}[!t]
\centering
\includegraphics[width=3.2in]{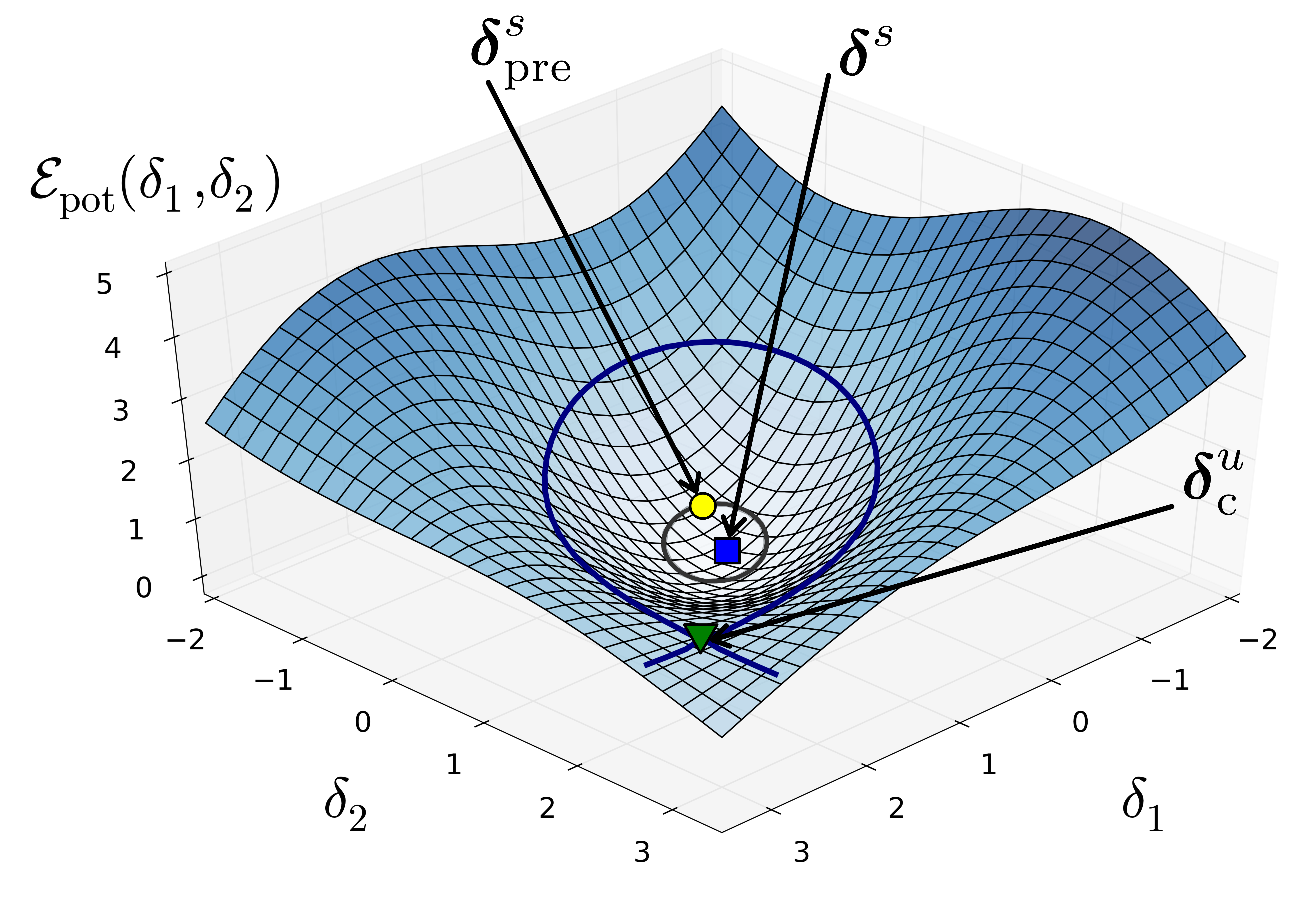}
\caption{(colour online) An illustration of the Hamiltonian (\ref{eq:ham}) for a TMIB system in the manifold where $\omega_1=\omega_2=0$. Equation (\ref{eq:potTMIB}) is plotted as a surface in $3$-dimensions and the energy difference $\delhh$ is given exactly by the difference in energy between the level sets $\mathcal{E}_\mathrm{pot}(\delta_1,\delta_2) = \mathcal{E}_\mathrm{c} = \mathcal{E}_\mathrm{pot}(\boldsymbol{\delta}_\mathrm{c}^u)$ and $\mathcal{E}_\mathrm{pot}(\delta_1,\delta_2) = \mathcal{E}_\mathrm{pot}(\boldsymbol{\delta}_\mathrm{pre}^s)$ where $\mathbf{x}_\mathrm{c}^u = [{\boldsymbol{\delta}_\mathrm{c}^u}^T,\mathbf{0}^T]^T$ is the closest UEP (found using condition (\ref{eq:cuep_def})), $\mathbf{x}_\mathrm{pre}^s = [{\boldsymbol{\delta}_\mathrm{pre}^s}^T,\mathbf{0}^T]^T$ is the pre-fault stable equilibrium point and $\mathbf{x}^s = [{\boldsymbol{\delta}^s}^T,\mathbf{0}^T]^T$ is the stable post-fault equilibrium point. \label{fig:cont1}}
\end{figure}

The expressions for kinetic and potential energy in the Hamiltonian function (\ref{eq:ham}) for this system are
\begin{equation*}
\begin{aligned}
&\mathcal{E}_\mathrm{kin}(\omega_1,\omega_2) = \frac{1}{2}M_1\omega_1^2 + \frac{1}{2}M_2\omega_2^2,
\end{aligned}
\label{eq:kinTMIB}
\end{equation*}
\begin{equation}
\begin{aligned}
&\mathcal{E}_\mathrm{pot}(\delta_1,\delta_2) = - (P_{\mathrm{m}1}-P_{\mathrm{a}1})\delta_1 - (P_{\mathrm{m}2}-P_{\mathrm{a}2})\delta_2  \\
&  \qquad - \bar{P}_{13}\cos(\delta_1) - \bar{P}_{23}\cos(\delta_2) - \bar{P}_{12}\cos(\delta_1-\delta_2).
\end{aligned}
\label{eq:potTMIB}
\end{equation}
where (\ref{eq:potTMIB}) is plotted as a surface in $3$-dimensions in Fig.~\ref{fig:cont1}. The critical energy boundary $\mathcal{E}_\mathrm{c}=\mathcal{H}(\mathbf{x}^u_\mathrm{c}) = \mathcal{E}_\mathrm{pot}(\boldsymbol{\delta}_\mathrm{c}^u)$ and the initial energy $\mathcal{H}(\mathbf{x}^s_\mathrm{pre})=\mathcal{E}_\mathrm{pot}(\boldsymbol{\delta}^s_\mathrm{pre})$ are plotted as level sets on the surface.

\section{An analytic stability metric}
\label{sec:low_bound}


%

An analytic stability metric, denoted as $\tau_\mathrm{A}$, which is purely a function of network parameters, is presented and derived here. The metric is formulated by considering (\ref{eq:solve_CCT}) and approximating both the Hamiltonian (\ref{eq:ham}) and the fault trajectory (\ref{eq:fault_traj}) in this equation as polynomial functions of rotor angles and time respectively. During a fault, it is assumed that the governor control systems for the mechanical input power $\Pmi$ in each generator are not able to act quickly enough to change the parameter value during (or immediately after) the fault; therefore $\Pmi^\mathrm{pre} = \Pmi^\mathrm{on} = \Pmi$ throughout this analysis. In addition, the dynamics of the rotor angles during the fault are approximated as constant but different accelerations. Although some of these approximations may seem cumbersome, and will detract significantly from the true dynamics of the system, they are valid in the limit as the CCT tends to zero. Therefore, we assume that for small values of the CCT these approximations can be assumed to capture the dynamics of a power system modelled as a Hamiltonian system. In addition, we remind the reader that this metric is designed to provide an instant illustration of the stability of a power system under the variation of a chosen parameter. In the next section, we demonstrate how our analytic metric can be used to find approximate regions for values of a load parameter in an aggregate network that improves stability for a given fault on the network. 

\subsection*{Metric formulation}
\label{subsec:metricform}
An analytical expression to approximate the CCT for three-phase to ground faults close to a given bus (on a balanced system such that it can be modelled by means of a single phase equivalent) can be found by adapting the energetic framework for the CCT presented in Section~\ref{sec:faults}. The expression (\ref{eq:solve_CCT}) is altered such that we solve
\begin{equation}
h_\mathrm{alt}(t) = \mathcal{E}_c,
\label{eq:new_ham_solve}
\end{equation}
where $h_\mathrm{alt}(t)$ is a polynomial function such that
\begin{equation}
h_\mathrm{alt}(t) \approx h(t),
\label{eq:key_ineq}
\end{equation}
with initial condition 
\begin{equation}
h_\mathrm{alt}(0) = h(0),
\label{eq:init_halt}
\end{equation}
and $h(t) \equiv \mathcal{H}(\mathbf{x}_\mathrm{on}(t))$.

In order to construct the function $h_\mathrm{alt}(t)$, we first approximate the Hamiltonian function as a polynomial function of the rotor angles, denoted $h_\mathrm{alt}(\boldsymbol{\delta}_\mathrm{on}(t))$. The kinetic term from the post-fault Hamiltonian function is removed by also modelling the dynamics during the fault as a Hamiltonian system. In general, there is no stable stationary point available during the fault so the power consumed by the conductive loads is approximated as a constant $P^\mathrm{on}_{\mathrm{a}i} = P^\mathrm{on}_{i}(\dspre)$ such that the dynamics can be written in the form (\ref{eq:genMMSEham}) to give
\begin{equation}
\dot{\mathbf{x}}_\mathrm{on} = \hat{\mathbf{F}}_\mathrm{on}(\mathbf{x}_\mathrm{on},\mathbf{x}^s_\mathrm{pre}),
\label{eq:ham_dyn}
\end{equation}
for $t\geq0$ and $\mathcal{H}(\mathbf{x}_\mathrm{on})\leq\mathcal{E}_c$. The Hamiltonian during the fault is given by
\begin{equation}
\mathcal{H}_\mathrm{on}(\mathbf{x}_\mathrm{on}(t)) =\mathcal{H}_\mathrm{on}(\mathbf{x}^s_\mathrm{pre}).
\end{equation}
In accordance with conservative systems, $\mathcal{H}_\mathrm{on}(\mathbf{x}^s_\mathrm{pre})$ is a constant in time and this property is used to recast the expression for $\mathcal{H}(\mathbf{x}_\mathrm{on}(t))$ in (\ref{eq:solve_CCT}) by considering the trivial relation 
\begin{equation}
\mathcal{H}(\mathbf{x}_\mathrm{on}(t)) = \mathcal{H}(\mathbf{x}_\mathrm{on}(t)) - \mathcal{H}_\mathrm{on}(\mathbf{x}_\mathrm{on}(t)) + \mathcal{H}_\mathrm{on}(\mathbf{x}^s_\mathrm{pre}),
\label{eq:Hamilts}
\end{equation}
resulting in the succinct expression
\begin{equation}
\begin{aligned}
\mathcal{H}(\mathbf{x}_\mathrm{on}(t)) & =
\sum_{i=1}^{n}(P_{\mathrm{a}i} - P_{\mathrm{a}i}^\mathrm{on})\delta_{\mathrm{on},i}(t) \\ & \quad - \sum_{\substack{i=1 \\ k>i}}^{n}(\bar{P}_{ik}-\bar{P}^\mathrm{on}_{ik}) \cos(\delta_{\mathrm{on},i}(t)-\delta_{\mathrm{on},k}(t)) \\ & \qquad + \mathcal{H}_\mathrm{on}(\mathbf{x}^s_\mathrm{pre}),
\end{aligned}
\label{eq:ham_fin}
\end{equation}
which has no dependence on rotor speeds. The values of the parameters $\bar{P}^\mathrm{on}_{ik}$, $P_{\mathrm{a}i}^\mathrm{on}$, $\bar{P}_{ik}$ and $P_{\mathrm{a}i}$ are found from the fault-on and post-fault reduced admittance matrices and the internal voltages.

A candidate function for $h_\mathrm{alt}(\boldsymbol{\delta}_\mathrm{on}(t))$ can be found by replacing the cosine terms in (\ref{eq:ham_fin}) with the function
\begin{equation*}
1-\frac{1}{2}\mathit{\Delta\delta}^2_{\mathrm{on},ik}(t) \approx \cos \left( \mathit{\Delta\delta}_{\mathrm{on},ik}(t)\right),
\end{equation*}
for small  $\mathit{\Delta\delta}_{\mathrm{on},ik}(t) = \delta_{\mathrm{on},i}(t) - \delta_{\mathrm{on},k}(t)$. This substitution gives
\begin{equation}
\begin{aligned}
& h_\mathrm{alt}(\boldsymbol{\delta}_\mathrm{on}(t)) = \sum_{i=1}^{n}(P_{\mathrm{a}i} - P_{\mathrm{a}i}^\mathrm{on})\delta_{\mathrm{on},i}(t) +  \\
&  \sum_{\substack{i=1 \\ k>i}}^{n}(\bar{P}_{ik}-\bar{P}^\mathrm{on}_{ik})
\left(
1-\frac{1}{2}\mathit{\Delta\delta}_{\mathrm{on},ik}^2(t)
\right) + \mathcal{H}_\mathrm{on}(\mathbf{x}^s_\mathrm{pre}) + C,
\end{aligned}
\label{eq:newadapt0}
\end{equation}
where the constant $C$ is found by applying the initial condition (\ref{eq:init_halt}), i.e. $h_\mathrm{alt}(\boldsymbol{\delta}_\mathrm{on}(0)) = \mathcal{H}(\mathbf{x}^s_\mathrm{pre})$ to (\ref{eq:newadapt0}). Therefore,
\begin{equation}
\begin{aligned}
& C = -\sum_{i=1}^{n}(P_{\mathrm{a}i} - P_{\mathrm{a}i}^\mathrm{on})\delta^s_{\mathrm{pre},i} +  \\
  & -\sum_{\substack{i=1 \\ k>i}}^{n}(\bar{P}_{ik}-\bar{P}^\mathrm{on}_{ik})
\left(
1-\frac{1}{2}\mathit{\Delta\delta}_{\mathrm{pre},ik}^2
\right) - \mathcal{H}_\mathrm{on}(\mathbf{x}^s_\mathrm{pre}) + \mathcal{H}(\mathbf{x}^s_\mathrm{pre}),
\end{aligned}
\label{eq:ham_const}
\end{equation}
where $\mathit{\Delta\delta}_{\mathrm{pre},ik}=\mathit{\Delta\delta}_{\mathrm{on},ik}(0)$ and (\ref{eq:newadapt0}) can be re-written as 
\begin{equation}
\begin{aligned}
& h_\mathrm{alt}(\boldsymbol{\delta}_\mathrm{on}(t)) = \sum_{i=1}^{n}(P_{\mathrm{a}i} - P_{\mathrm{a}i}^\mathrm{on})(\delta_{\mathrm{on},i}(t)-\delta^s_{\mathrm{pre},i}) +  \\
&  \sum_{\substack{i=1 \\ k>i}}^{n}\frac{(\bar{P}_{ik}-\bar{P}^\mathrm{on}_{ik})}{2}
\left(
\mathit{\Delta\delta}_{\mathrm{on},ik}^2(t)-\mathit{\Delta\delta}_{\mathrm{pre},ik}^2
\right) + \mathcal{H}(\mathbf{x}^s_\mathrm{pre}),
\end{aligned}
\label{eq:newadapt}
\end{equation}
where the constant $C$ is written explicitly.

In order to make (\ref{eq:newadapt}) an explicit polynomial function of time, the fault trajectory must also be written as a polynomial function of time. In general, the dynamics during a fault are non-trivial \cite{Chu1997} but in order for (\ref{eq:new_ham_solve}) to be analytically solvable for time, the rotor angle dynamics in (\ref{eq:newadapt}) must have the form
\begin{equation}
\delta_{\mathrm{on},i}(t) = \frac{1}{2}u_i t^2 + \delta^s_{\mathrm{pre},i},
\label{eq:diff_acc}
\end{equation}
where the initial condition $\dot{\delta}_{\mathrm{on},i}(0)=0$ holds for all $i$. An appropriate value for the acceleration $u_i$ can be found by assuming that for small CCTs the rotor dynamics can be modelled as a constant acceleration equal to the initial rotor acceleration at $t=0$. This is given by 
\begin{equation}
\dot{\mathbf{x}}_\mathrm{on} \approx \hat{\mathbf{F}}_\mathrm{on}(\mathbf{x}^s_\mathrm{pre},\mathbf{x}^s_\mathrm{pre}) 
= 
\begin{pmatrix}
\boldsymbol{\omega} \\
\hat{\mathbf{A}}_\mathrm{on}(\dspre,\dspre)
\end{pmatrix}
=
\begin{pmatrix}
\boldsymbol{\omega} \\
\mathbf{A}_\mathrm{on}(\dspre)
\end{pmatrix},
\label{eq:vectondyn}
\end{equation}
for short fault times. From equation (\ref{eq:vectondyn}) the rotor accelerations $\ddot{\delta}_{\mathrm{on},i} = u_i =   A_{\mathrm{on},i}(\dspre)$. By substituting expressions (\ref{eq:diff_acc}) for the rotor angles into (\ref{eq:newadapt}), the function
\begin{equation}
\begin{aligned}
h_\mathrm{alt}(t) &= \mathcal{H}(\mathbf{x}^s_\mathrm{pre}) +\sum_{i=1}^{n}(P_{\mathrm{a}i} - P_{\mathrm{a}i}^\mathrm{on})\frac{1}{2}u_i t^2 + \\ 
& \sum_{\substack{i=1 \\ k>i}}^{n}(\bar{P}_{ik}-\bar{P}^\mathrm{on}_{ik})
\left(
\frac{1}{8}u_{ik}^2t^4 + \frac{1}{2}u_{ik}\mathit{\Delta\delta}_{\mathrm{pre},ik}t^2\right),
\end{aligned}
\label{eq:HaFin}
\end{equation}
is a quadratic in $t^2$ where $u_{ik} = u_i - u_k$. Now (\ref{eq:new_ham_solve}) can be written as
\begin{equation} 
\alpha t^4 + \beta t^2 - \gamma = 0,
\label{eq:quadratic1}
\end{equation}
where the coefficients 
\begin{eqnarray*}
\alpha &=& \sum_{\substack{i=1 \\ k>i}}^{n}
\frac{1}{8}(\bar{P}_{ik}-\bar{P}^\mathrm{on}_{ik})u_{ik}^2, \\ 
\beta &=&  \sum_{\substack{i=1 \\ k>i}}^{n}{\frac{1}{2}(\bar{P}_{ik}-\bar{P}^\mathrm{on}_{ik})u_{ik}\mathit{\Delta\delta}_{\mathrm{pre},ik}} + \\
& & \qquad \sum_{i=1}^{n}\frac{1}{2}(P_{\mathrm{a}i} - P_{\mathrm{a}i}^\mathrm{on})u_i, \\ 
\gamma &=&  \mathcal{E}_\mathrm{c} - \mathcal{H}(\mathbf{x}_\mathrm{pre}^s) = \delhh > 0
\end{eqnarray*}
are functions of the power network parameters. The solution of (\ref{eq:quadratic1}) and thus the expression for our analytic CCT is given by
\begin{equation}
\tau_\mathrm{A} = \left( \frac{ -\beta \pm \sqrt{\beta^2 + 4\alpha\gamma }}{2 \alpha} \right)^\frac{1}{2}.
\label{eq:final_res}
\end{equation}
The smallest real value of $\tau_\mathrm{A}$ is taken for a given set of parameters. A purely imaginary value for $\tau_\mathrm{A}$ is produced if the discriminant $\beta^2 + 4\alpha\gamma <0$ or if $\beta<0$ and $\alpha <0$. In the case where $\alpha <0$ and $\beta >0$ two positive roots are produced, otherwise there is one real root to (\ref{eq:final_res}). However, in general the parameter $\alpha$ is positive because the total electrical load of a network reduces during a fault and so it is reasonable to assume that $\bar{P}_{ik} > \bar{P}^\mathrm{on}_{ik}\; \mathrm{for\;all}\;i,k$.

\begin{figure}[!t]
\centering
\includegraphics[width=3.4in]{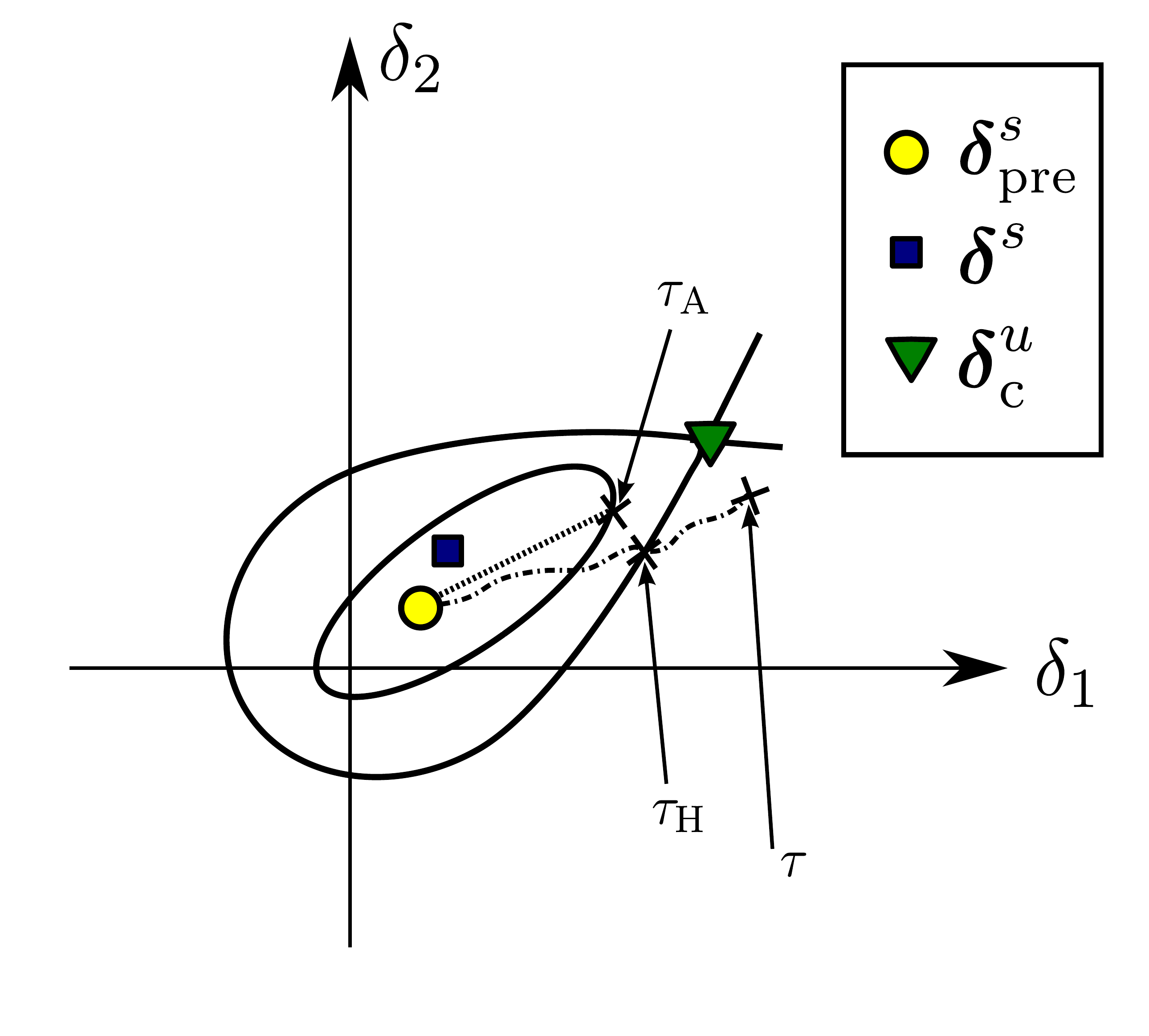}
\caption{(colour online) This figure illustrates the definitions of the three CCT metrics in $\delta$-space: $(i)$ the true CCT $\tau$ is the maximum time a fault can remain on-line such that there is one full swing of the rotor angles post-fault; $(ii)$ the CCT estimate $\tau_\mathrm{H}$ is defined using the direct methods by solving (\ref{eq:solve_CCT}) and in the figure it is where the fault trajectory (dashed line) intersects the level set $\mathcal{H}(\mathbf{x})=\mathcal{E}_c$; $(iii)$ the analytic CCT $\tau_\mathrm{A}$ is the analytic metric derived in this paper. It is the solution to (\ref{eq:new_ham_solve}), where the fault-on and post-fault regimes are modelled as hamiltonian systems. In $\delta$-space (\ref{eq:new_ham_solve}) is given by an ellipse for a TMIB system and the fault trajectory (dotted line) is approximated using a constant but unique acceleration for each generator. \label{fig:tcct}}
\end{figure}

Figure~\ref{fig:tcct} illustrates how the analytic CCT $\tau_\mathrm{A}$ compares with the true CCT $\tau$ and the CCT estimate $\tau_\mathrm{A}$ developed in Section~\ref{sec:faults}.

\section{Parametric stability analysis}
\label{sec:3genExample}

\subsection{Implementation details}

The stability of a power network is not only dependent on the type or duration of a fault but also on the choice of system parameters. Optimal regions of parameter space that increase the stability of a power system can be identified by the variation of system parameters. Here, we investigate values for a load on a TMIB network which improve system stability for a given fault, by comparing the metrics $\tau_\mathrm{H}$ and $\tau_\mathrm{A}$ outlined in Sections~\ref{sec:fault} and \ref{sec:low_bound} against the true CCT $\tau$. The energy difference $\delhh$ is also compared against the temporal stability metrics. 

The variation of a load parameter in a power network will change one of the elements in the full bus admittance matrix $\mathbf{Y}_\mathrm{BUS}$, but due to Kron reduction all the parameters in the reduced admittance matrix $\mathbf{Y}_\mathrm{red}$ in each regime will change. Therefore, for each incremental change in the parameter we consider, a new fault study is required to find the CCT. In each fault study, the mechanical input powers for each generator are found by performing a pre-fault power flow for the system. The power flow is conducted using the same (small) rotor angles found in \cite{Anderson2002book} for each incremental change in the parameter value to ensure that the network is initially in a stable state. 

It is relatively quick to conduct a single fault study to find the true CCT $\tau$. However, for each incremental change in a parameter, a new fault analysis is required to find the CCT. An analytic CCT $\tau_\mathrm{A}$ has been developed to study trends in stability of power systems, which can be found instantly once the relevant parameter values for the model in (\ref{eq:fault_model}) have been collected. For a given fault, system parameters can be varied continuously using an analytic CCT and this can provide an initial picture of stability that informs more detailed analysis.

All the stability metrics introduced in this paper, except for the true CCT $\tau$, are dependent on a critical energy boundary $\mathcal{E}_\mathrm{c}$ which is dependent on the location of the closest UEP in this work. The position of the closest UEP will change under the variation of the loads and there are techniques developed in the literature to find these quickly \cite{Liu1997}. However, we choose to use numerical continuation (previously applied to power systems in \cite{Chen2009}) to illustrate interesting features of the closest UEP under the variation of loads. The stationary points of a TMIB system, modelled by the ODEs in (\ref{eq:TMIBode}), are located using the continuation software AUTO\footnote{http://indy.cs.concordia.ca/auto/} as a load parameter is varied. There is no rigorous proof provided in this paper that all the possible unstable equilibria on the stability boundary of a stable equilibrium point can be found from the solution branches from numerical continuation. However, no other solutions were found for this system when performing an exhaustive search over state space using the root finding algorithm $\mathrm{fsolve}$ from the Scipy\footnote{http://docs.scipy.org/doc/} library in the Python\footnote{www.python.org} programming language. Therefore, without further analysis, it is assumed that only in a TMIB system can all the necessary stationary points be found using this continuation method.




The stationary points for each value of the continuation parameters in Fig.~\ref{fig:CCTbifFOU}a and Fig.~\ref{fig:CCTbifFOUi}a are obtained by the following method:  A stable stationary point denoted by ${\mathbf{x}^s = [\delta_1^s, \delta_2^s,\omega_1=0,\omega_2=0]}$ that solves the post fault equation $\mathbf{F}_\mathrm{post}(\mathbf{x}^s) = \mathbf{0}$ (where $n=2$ and $\delta_3=0$) is found using the root finding algorithm $\mathrm{fsolve}$, where $|\delta_i^s-\delta_k^s| \in \frac{\pi}{2}\,\mathrm{for}\,i=1,2$. This point belongs to the lower branch (blue squares) of the bifurcation diagram in Figs~\ref{fig:CCTbifFOU}a and \ref{fig:CCTbifFOUi}a.  The other (unstable) equilibria on the boundary of the stability region of the stable equilibrium point $\mathbf{x}^s$ which satisfy $\hat{\mathbf{F}}_\mathrm{post}(\mathbf{x})=\mathbf{0}$ are found by numerical continuation of the element $B_{12}$ from the reduced admittance matrix $\mathbf{Y}_\mathrm{red}$. Once the continuation branches are found, the stationary points at the value of $B_{12}$ found in $\mathbf{Y}_\mathrm{red}$ are recorded. The local stability of the stationary points obtained are found by computing the eigenvalues of the Jacobian matrix for the system (\ref{eq:TMIBode}). The stability of the solution branches in Figs~\ref{fig:CCTbifFOU}a and \ref{fig:CCTbifFOUi}a are stated in terms of the number of eigenvalues with real part greater than zero which can be found in the figure caption of Fig~\ref{fig:CCTbifFOU}.

\begin{figure}[!t]
\centering
\includegraphics[width=3.4in]{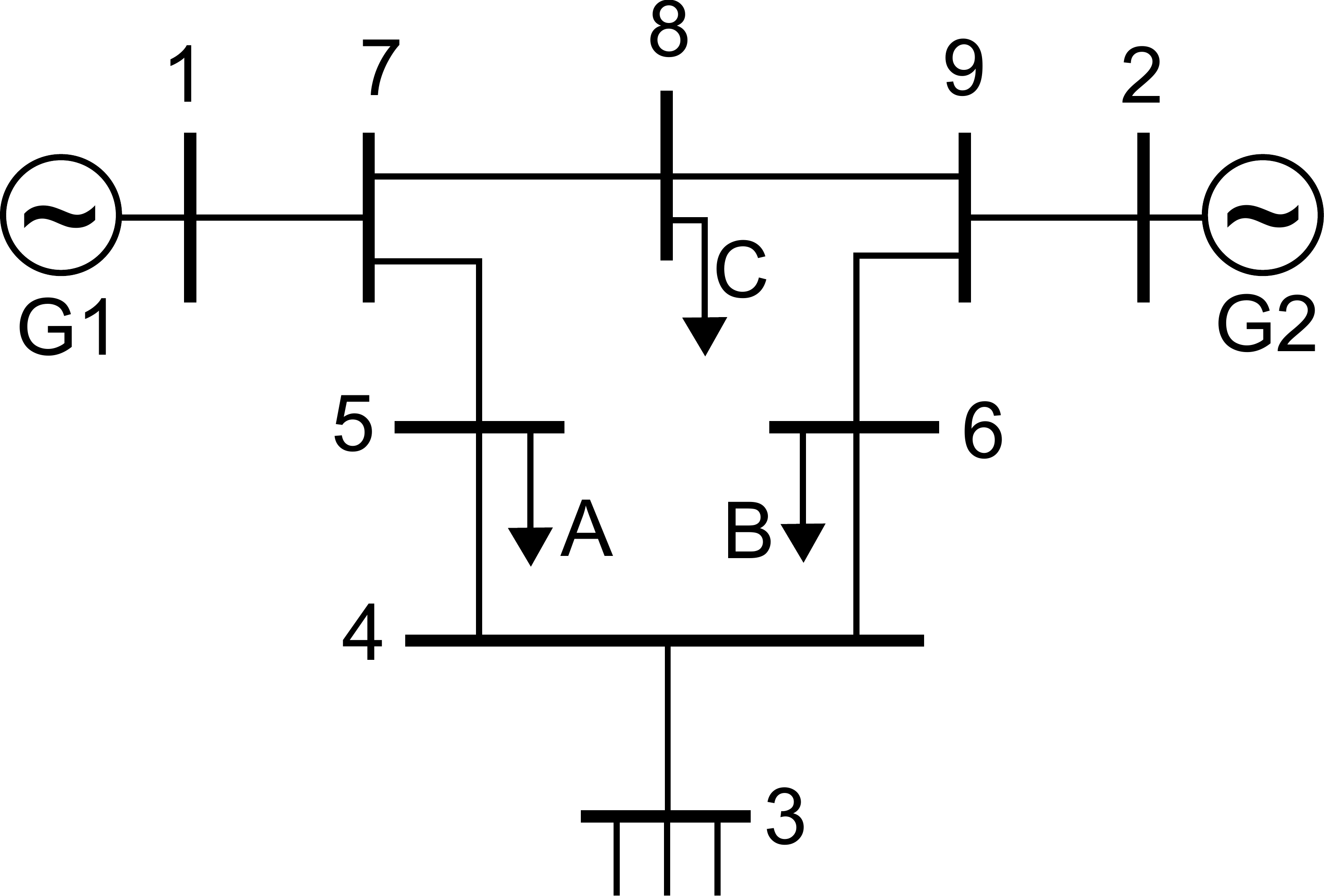}
\caption{Schematic of the full test network found in \cite{Anderson2002book} where all parameter values can be found in the reference. Buses 1 and 2 are $(P,V)$ buses with synchronous generators attached. Bus 1 is a $(V,\theta)$ infinite bus. All other buses are $(P,Q)$ buses with buses 5, 6 and 8 also possessing shunt loads with the original values $Y_A = 1.261 -0.2634j$, $Y_B = 0.8777 - 0.0346j$ and $Y_C = 0.969-0.1601j$. The fault we consider occurs on the line 5-7 close to bus 7, and the post fault network has line 5-7 switched out. \label{fig:FouadNetwork}}
\end{figure}

In this study a 9 bus, 3 generator power network found in \cite{Anderson2002book} is used, a schematic of this network is provided in Fig.~\ref{fig:FouadNetwork}. All parameter values for this network are taken from \cite{Anderson2002book}. This network is adapted into a TMIB network by changing one of the generators, which has an inertia an order of magnitude larger than the other two, into an infinite bus. The specific fault we consider is a three-phase to ground fault close to bus 7 on the line connecting buses 5 and 7. The post fault network is identical to the pre-fault network except the line connecting buses 5 and 7 is switched out.

\subsection{Results}

\begin{figure}[!t]
\centering
\includegraphics[width=3.4in]{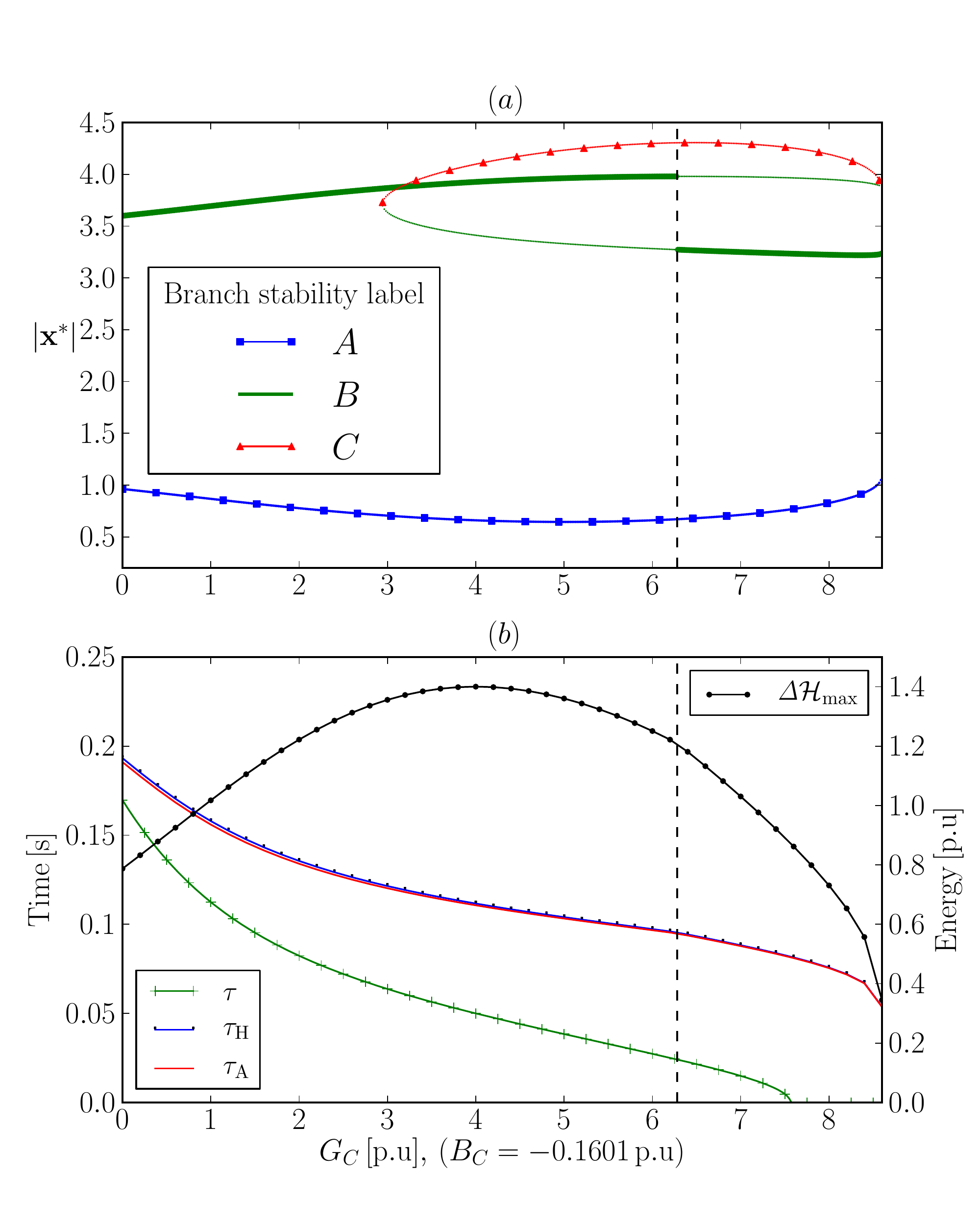}
\caption{(colour online) The continuation diagram in (a) plots the modulus of the stationary solutions $|\mathbf{x}^*|$ to the ODE (\ref{eq:TMIBode}) as a function of the bifurcation parameter $G_C$ with $B_C=-0.1601$ p.u. The stability of the solution branches are colour coded using the legend. The stability of each branch is given by the number of eigenvalues with positive real part; for branch segments $A$, $B$ and $C$ these are $0$, $1$ and $2$ respectively. The thicker line indicates the closest UEP $|\mathbf{x}^u_\mathrm{c}| = |\boldsymbol{\delta}^u_\mathrm{c}|$. In (b) the CCT metrics $\tau$, $\tau_\mathrm{H}$ and $\tau_\mathrm{A}$ should be read using the left y-axis and the energy margin metric $\delhh$ should be read using the right y-axis. Note that the values of $\tau_\mathrm{H}$ and $\tau_\mathrm{A}$ are very close together. \label{fig:CCTbifFOU}}
\end{figure}

\begin{figure}[!t]
\centering
\includegraphics[width=3.4in]{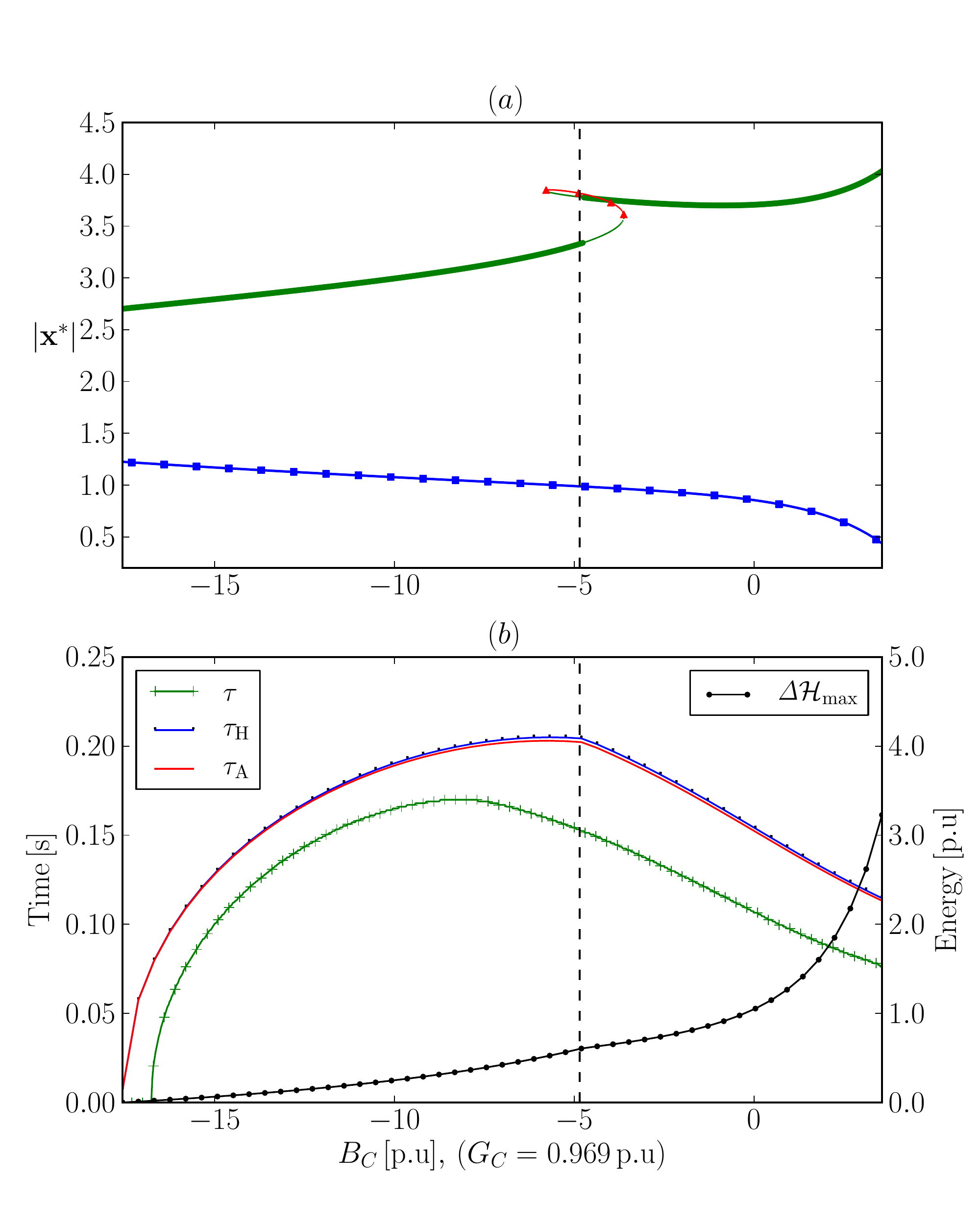}
\caption{(colour online) The continuation diagram in (a) plots the modulus of the stationary solutions $|\mathbf{x}^*|$ to the ODE (\ref{eq:TMIBode}) as a function of the bifurcation parameter $B_C$ with the conductance of load C held constant at $G_C=0.969$ p.u. The stability of the solution branches are colour coded using the legend in Fig.~\ref{fig:CCTbifFOU}a and the stability information can be found in the caption of Fig.~\ref{fig:CCTbifFOU}. The thicker line indicates the closest UEP $|\mathbf{x}^u_\mathrm{c}| = |\boldsymbol{\delta}^u_\mathrm{c}|$. In (b) the CCT metrics $\tau$, $\tau_\mathrm{H}$ and $\tau_\mathrm{A}$ should be read using the left y-axis and the energy margin metric $\delhh$ should be read using the right y-axis. Note that the values of $\tau_\mathrm{H}$ and $\tau_\mathrm{A}$ are very close together. \label{fig:CCTbifFOUi}}
\end{figure}

The size and nature of loads on actual power systems can vary over time and, in respect of the susceptive part, can be modified by the addition of reactive compensation. As a consequence, the conductive and susceptive parts of load C (denoted $G_C$ and $B_C$ respectively) of the network in Fig.~\ref{fig:FouadNetwork} are investigated by varying one part while maintaining the other constant at its original value. In Figs.~\ref{fig:CCTbifFOU}a and \ref{fig:CCTbifFOUi}a the domains of the parameters $B_{C}$ and $G_{C}$ respectively are constrained by two conditions: $(i)$ the energy margin $\delhh \geq 0$ and $(ii)$ that the synchronous machines are operating as generators in the pre-fault power flow, i.e. $\Pmo>0$ and $\Pmt>0$. There is an additional constraint in Fig.~\ref{fig:CCTbifFOU} where only positive values of conductance are explored. 

The critical energy change for the system $\delhh$ (black line), the CCT estimate $\tau_\mathrm{H}$ (blue line) and the analytic CCT $\tau_\mathrm{A}$ (red line) are plotted as functions of the continuation parameters in the lower panels of Figs.~\ref{fig:CCTbifFOU} and \ref{fig:CCTbifFOUi}. In addition, the true CCT $\tau$ (green line) is plotted using a simple algorithm that uses a binary search to find the maximum duration which the fault can be left on-line such that the rotor angles have one full swing together before they diverge. There are two different scales to facilitate observing the functions in the lower panels of Figs.~\ref{fig:CCTbifFOU} and \ref{fig:CCTbifFOUi}: the energy change  $\delhh$ should be read using the right-hand y-axis labels and the three time metrics should be read using the left-hand y-axis labels as indicated in the figures.

\begin{figure}[!t]
\centering
\includegraphics[width=3.4in]{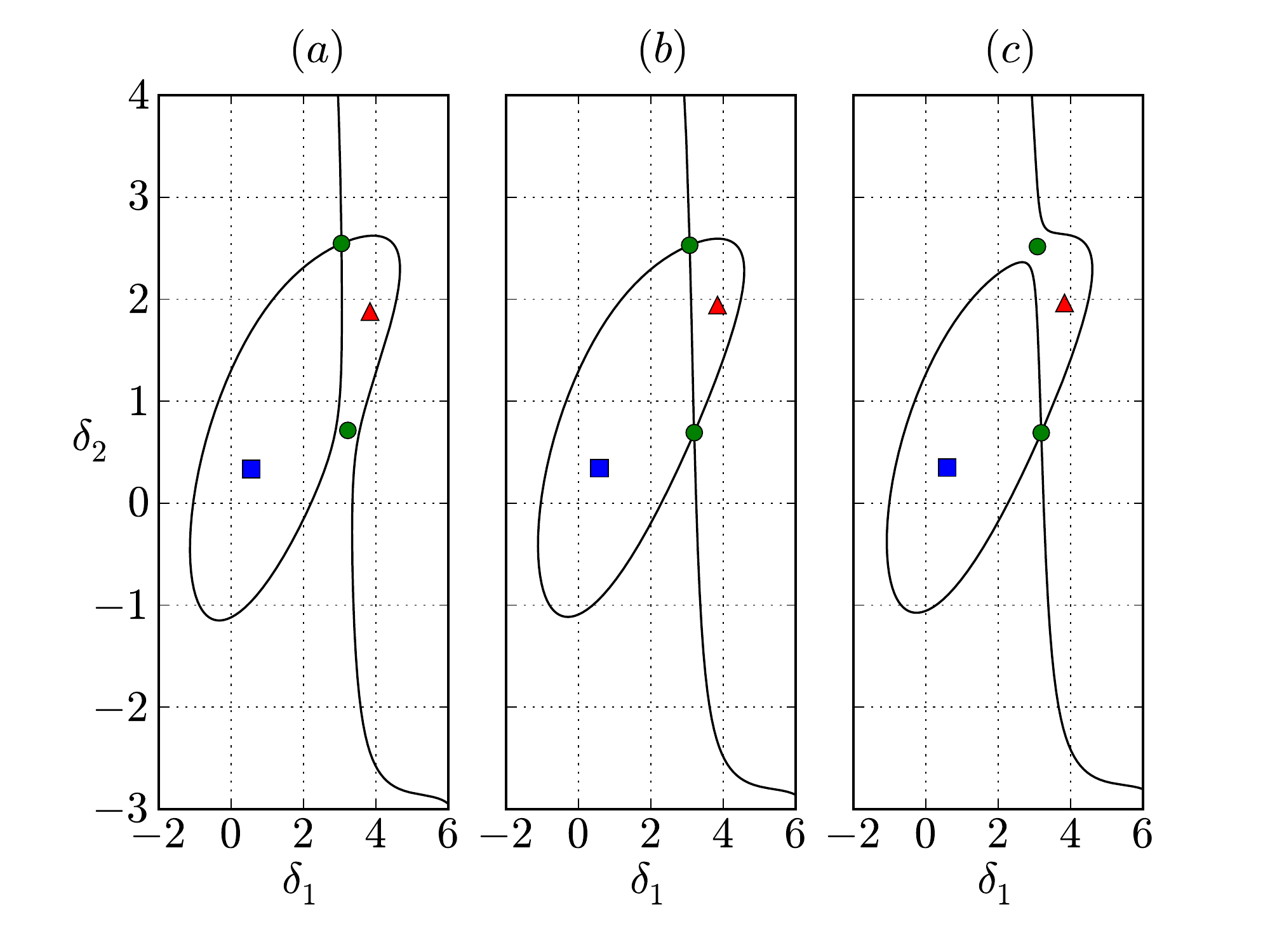}
\caption{(colour online) This figure illustrates the change to the energy boundary as the closest UEP changes position at the discontinuity $\hat{G}_\mathrm{C}= 6.28$ (vertical dotted line) in Fig.~\ref{fig:CCTbifFOU}a due to a change in the parameter $G_\mathrm{C}$. The energy boundary is plotted in the manifold where $\omega_1=\omega_2=0$ such that the boundary can be plotted by the level sets $\mathcal{E}_\mathrm{pot}(\delta_1,\delta_2) = \mathcal{E}_\mathrm{c}$. These level sets are plotted for conductance values (a) $G_\mathrm{C} = 5.5 < \hat{G}_\mathrm{C}$, (b) $G_\mathrm{C}= \hat{G}_\mathrm{C}$ and (c) $G_\mathrm{C} = 6.6 > \hat{G}_\mathrm{C}$. In each sub-figure, the stationary points are plotted using the same marker style found in the legend of Fig.~\ref{fig:CCTbifFOU}a.\label{fig:splices}}
\end{figure}

In Fig.~\ref{fig:CCTbifFOU} the dependence of the system stability, for the fault we consider, on the conductance $G_C$ is studied as the susceptance $B_C$ is held constant. In Fig.~\ref{fig:CCTbifFOU}a there is a discontinuity in the closest UEP (thick line) at $G_C = 6.26$ (vertical dotted line) which is located between two pairs of fold points at $G_C = 2.95$ and $G_C = 8.56$. (See \cite{strogatz2001nonlinear} for an explanation of fold points). In Fig.~\ref{fig:CCTbifFOU}b there is a discontinuous change in the gradient of $\delhh$ which coincides with the discontinuity at $G_{C} = 6.26$, but the maximum point for $\delhh$ at $G_{C} =4.0$  does not coincide with the other discontinuity in the closest UEP, nor any other points of significance in Fig.~\ref{fig:CCTbifFOU}a. The analytic CCT $\tau_\mathrm{A}$ is observed to approximate the CCT estimate $\tau_\mathrm{H}$ very well as the load parameter $G_C$ is varied. However, the energetic techniques used to find $\tau_\mathrm{A}$ and $\tau_\mathrm{H}$ have resulted in non-conservative estimates of the true CCT $\tau$. 
This feature is a manifestation of the original issue with direct methods which concerns the dissipative term (\ref{eq:cond_elec}) at each bus in the reduced network. Direct methods can be used for conservative stability assessments where the transfer conductances $G_{ij}$ in the reduced network matrix $\mathbf{Y}_\mathrm{red}$ are assumed to be small or zero \cite{Chiang2011book}, therefore the CCT estimate $\tau_\mathrm{H}$ and the analytic CCT $\tau_\mathrm{A}$ are strict lower bounds of the true CCT $\tau$ for networks with zero transfer conductances. However, even for a network with lossless lines, Kron reduction invokes complications in which a shunt load conductance in the full bus admittance matrix will increase the absolute values of $G_{ij}$ in the reduced admittance matrix \cite{Ribbens-Pavella1985, Sastry1973}. 

\begin{figure}[h]
\centering
\includegraphics[width=3.4in]{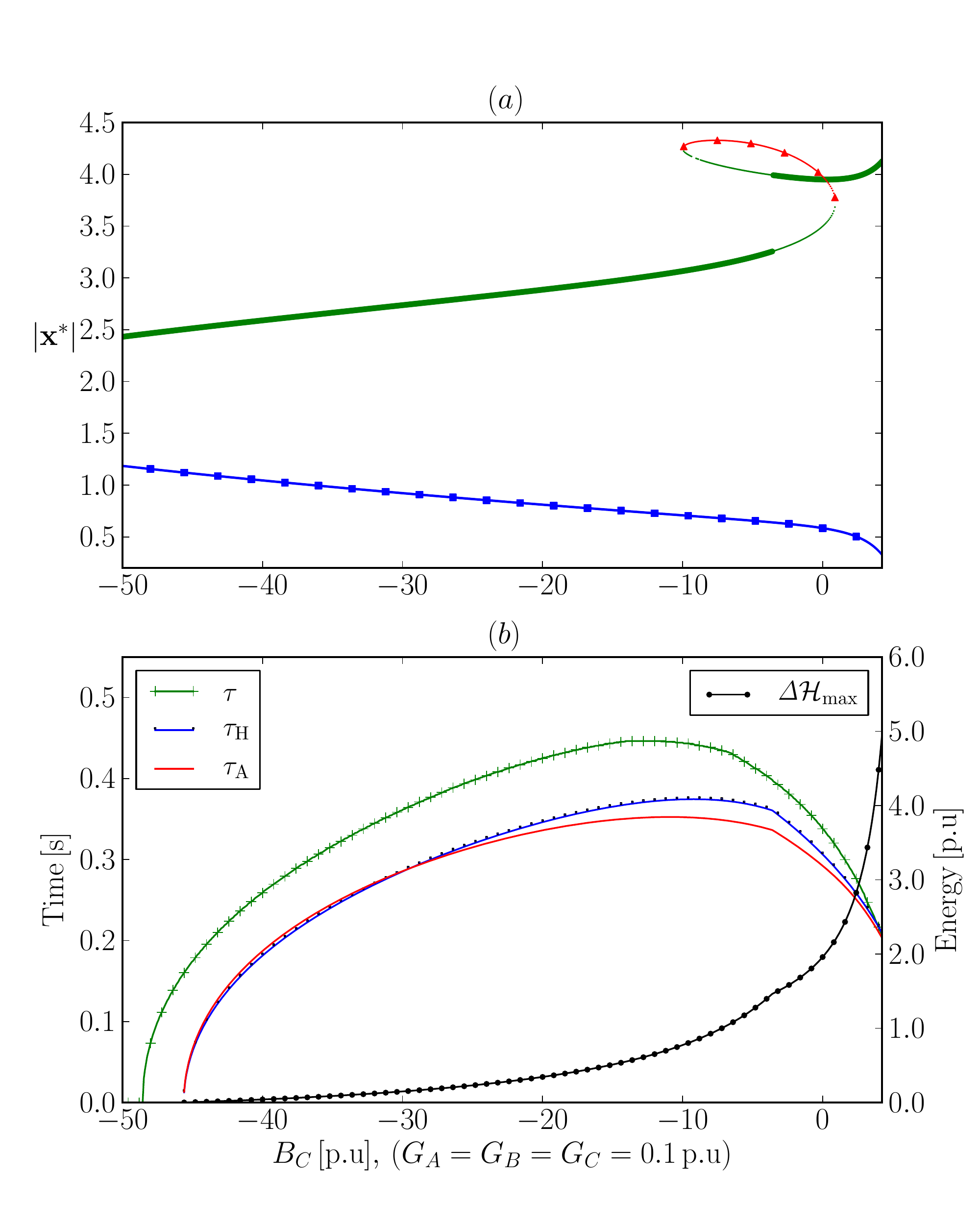}
\caption{(colour online) This Figure is identical to the results presented in \ref{fig:CCTbifFOUi} except that the conductive parts of the loads in the network are changed to be small values. This is done to show that the analytical CCT serves a lower bound for networks with small transfer conductances in the reduced admittance matrix and this is achieved by having small conductive parts to the loads. The values are $G_A =G_B=G_C= 0.1$ for loads A, B and C respectively. 
\label{fig:CCTbifFOUiSmall}}
\end{figure}

In Fig.~\ref{fig:CCTbifFOUiSmall}, the susceptance $B_C$ is varied in a network with small load conductances and it is observed that $\tau_\mathrm{A}$ and $\tau_\mathrm{H}$ are lower bounds to $\tau$ when compared to the results in Fig.~\ref{fig:CCTbifFOUi}. Despite whether the analytic CCT is an over or underestimate of the true CCT, it performs well as an indicator of the expected increase or decrease of the CCT as $G_C$ is changed. The greatest CCT as measured by all metrics for the fault we have considered is produced at $G_C=0$. (This behaviour was found for all possible faults on the network under the variation of one of the loads $A$, $B$ or $C$ within an order of magnitude of its nominal value.) In general, a lower mechanical input power from each generator is required for lower network loadings and therefore the acceleration of the generator rotor angles is roughly proportional to the mechanical input power, assuming that the load of the network decreases during a fault. Therefore, there is more time for the rotors to reach a critical value where they begin to diverge. 

In Fig.~\ref{fig:splices} the change to the energy boundary $\mathcal{H}(\mathbf{x})=\mathcal{E}_\mathrm{c}$ (plotted in the manifold where $\omega_1=\omega_2=0$) due to the discontinuity in the location of the closest UEP is illustrated in Fig.~\ref{fig:CCTbifFOU}a. In this manifold, the energy boundary is plotted as the level set $\mathcal{E}_\mathrm{pot}(\delta_1,\delta_2) = \mathcal{E}_\mathrm{c}$ (black line) for conductance values $(a)$ $G_C= 5.5 < \hat{G}_C$, $(b)$ $G_C= \hat{G}_C$ and $(c)$ $G_C= 6.6 > \hat{G}_C$ where the discontinuity occurs at $\hat{G}_C = 6.28$. In each sub-figure of Fig.~\ref{fig:splices} the stationary points of (\ref{eq:TMIBode}) are plotted using the same marker style found in the legend of Fig.~\ref{fig:CCTbifFOU}a and the level set is observed to intersect the closest UEP.

In Fig.~\ref{fig:CCTbifFOUi} the dependence of the system stability, for the fault we consider, on the susceptance $B_C$ is studied as the conductance $G_C$ is held constant. In Fig.~\ref{fig:CCTbifFOUi}a there is a discontinuity in the closest UEP (thick line) at $B_C = - 4.80$ which is located between two fold points at $B_C = -5.78$ and $B_C = -3.62$. In Fig.~\ref{fig:CCTbifFOUi}b there is a discontinuous change in the gradient of $\delhh$, $\tau_\mathrm{H}$ and $\tau_\mathrm{A}$ at the discontinuity in the closest UEP. The maximum of $\delhh$ occurs at the highest value of susceptance plotted, which shows that the energy margin is not the best metric to quantify stability. The analytic CCT is very close to the CCT estimate as $B_C$ is varied and are, again, overestimates due to the presence of non-negligible transfer conductances. The maximum points of $\tau_\mathrm{H}$ and $\tau_\mathrm{A}$, both at $G_C = -5.75$ are very close to the closest UEP discontinuity, however the maximum point of the true CCT $\tau$ is lower, at $B_C = -8.2$. Despite this, the change in the true CCT as the susceptance $B_C$ is varied is well captured by the CCT approximations, except in the region $[-8.2,-4.8]$ where the gradients are of different signs. 

\begin{table}[h!]
\begin{center}
\begin{tabular}{llll}
\hline
 & & \\
Optimum susceptance   &  $\tau(\cdot)$  &   $\tau_\mathrm{A}(\cdot)$  \\
\hline 
 & & \\
    \vspace{2mm} 
 $\argmax (\tau,B_B) = -3.10$   &  $0.124$s & N/A \\
      \vspace{2mm} 
 $\argmax (\tau_\mathrm{A},B_B)=-1.20$   &  $0.120$s & $0.157$s \\
 \hline
\end{tabular}
\end{center}
\caption{Maximum CCT values at optimum susceptance values for {load B}. The original true CCT at the parameter values given in \cite{Anderson2002book} is $\tau = 0.107s$}
\label{tab:optLoadB}
\end{table}
\begin{table}[h]
\begin{center}
\begin{tabular}{lll}
\hline
 & & \\
Optimum susceptance   &  $\tau(\cdot)$  &   $\tau_\mathrm{A}(\cdot)$  \\
\hline 
  & & \\
 \vspace{2mm}
  $\argmax (\tau,B_C)= -8.20$   &  $0.170$s & N/A \\
      \vspace{2mm}
  $\argmax (\tau_\mathrm{A},B_C)=-5.75$   & $0.158$s  & $0.203$s \\
\hline
\end{tabular}
\end{center}
\caption{Maximum CCT values at optimum susceptance values for {load C}. The original true CCT at the parameter values given in \cite{Anderson2002book} is $\tau = 0.107s$}
\label{tab:optLoadC}
\end{table}
It is observed that the system stability can benefit by setting the susceptance of load C to the optimum susceptance as measured by the analytic CCT because these susceptance values can be evaluated without the need for numerical integration. From Table~\ref{tab:optLoadC}, the true CCT using the parameter values as stated in \cite{Anderson2002book} is $\tau = 0.107\mathrm{s}$. A network operating at the optimum value of susceptance $B_C = -8.2$ would give a true CCT of $\tau(B_C = -8.2) = 0.170\mathrm{s}$ and this is an increase of $0.63\mathrm{s}$. However, the value of the true CCT at the optimum value of susceptance as measured by the analytic CCT is $\tau(B_C = - 5.75) = 0.158\mathrm{s}$ which is a smaller but significant increase of $0.51\mathrm{s}$. The advantage of using the optimum values of susceptance as measured by the analytic CCT is that an improved susceptance value is known as soon as the relevant network parameters have been collected.

The continuation of the susceptive part of load B $B_B$ (with the other loads at original values) is considered for the same fault at bus $7$ and the results are given in Table~\ref{tab:optLoadB}. The results for load $A$ are not included in the tables because there was no maximum point for CCT found as the susceptive part of load A $B_A$ was varied and the trends in stability were similar to the lower panel of Fig~\ref{fig:CCTbifFOU}. The largest CCT was found for $B_A = -13.9$ which is the lowest value of susceptance for which the mechanical input powers of the synchronous generators $\Pmo$ and $\Pmt$ are both greater than zero.



\section{Discussion}
\label{sec:discussion}

In this paper we have presented a new analytic CCT metric $\tau_\mathrm{A}$ designed to be able to capture trends in the true CCT $\tau$ as a system parameter is varied. Specifically the effects of the conductive and susceptive parts of a load parameter on the network were considered, given a fault on the network. The analytic CCT metric $\tau_\mathrm{A}$ was formulated by taking a polynomial approximation of the CCT estimate $\tau_\mathrm{H}$ developed from direct methods and it is found that despite the simplicity of the formula, the analytic CCT is a good approximation to $\tau_\mathrm{H}$ for short times. Given the difficulty for all direct methods to incorporate power networks with non-negligible transfer conductances, it was expected that the two approximating metrics $\tau_\mathrm{H}$ and $\tau_\mathrm{A}$ would not perform well as good estimates to the true CCT $\tau$. However, the two approximating metrics performed much better as indicators of trends in stability as a load parameter was varied.  In addition, the results in this paper were generated for an aggregate power network, the two-machine infinite-bus, where studies on aggregated networks \cite{Johnstone2014,Hughes2006} are generally conducted to analyse global trends of a network with a much larger number of generators, buses and loads. The approximating CCT metrics are valid in principle for power systems with or without an infinite bus and numerical results will be extended to alternative networks without an infinite bus in future work. 

Direct methods were chosen to formulate the CCT approximations because of their ability to construct a well-defined stability boundary in terms of a critical energy $\mathcal{E}_\mathrm{c}$. However, this stability boundary is dependent on finding a critical UEP of the system that can be used to approximate the energy when a power system is modelled as a Hamiltonian system.  In this paper, the closest UEP was chosen to compute the energy boundary because it is valid for any fault on a power system. A more accurate method to quantify the system energy uses the controlling UEP \cite{Chiang1989a}, which is dependent on the fault a network suffers and can be found using the `boundary of stability region-based controlling unstable equilibrium point' (BCU) method \cite{Chiang2011book,Chiang1994}; the limitations of which are discussed in \cite{Llamas1995}. For the TMIB network it was possible to find all UEPs on the stability boundary of an appropriate stable equilibrium point of the system under the variation of a load parameter using numerical continuation methods. There was no attempt to formally prove that all UEPs are captured, but an exhaustive algorithmic search was used to confirm this. However, the scalability of this approach is limited because the number of equilibria increases as the system size increases \cite{Tavora} and it is increasingly difficult to identify the critical UEP and this is a significant area of research in itself \cite{Liu1997,Lee2003}. In addition, the positions of the equilibria have to be found for each incremental change of a load parameter. This is another reason the analysis in this paper was limited to study trends on an aggregated power network with a small number of generators.

The most general use for the analytic CCT proposed in this paper is to capture stability trends under the variation of a network or generator parameter and we have specifically studied the effect of varying loads. A more specific suggested use for the metric is to locate regions of parameter space that will improve the system stability in terms of CCT. For the fault studied in our analysis of the TMIB network, we found that the optimum value of CCT for the conductive part of load C is $G_C=0$ for non-negative conductance values and that the CCT decreases as $G_C$ increases. One interpretation of this result is that a high penetration of linear power sources, local to the point of power consumption has the effect of increasing system stability, but research in \cite{Slootweg2002} suggests that the effect on system dynamics is dependent on the specific technologies in the generator. More promising results were found for the variation of the susceptive part of a load under constant conductance. The variation of susceptive loads can represent, for example, a network owner's installation of reactive compensation, a measure that is known to contribute not only to voltage regulation but also transient stability \cite{Delfino1988}. An optimal value of susceptance that maximises the CCT was identified in all three temporal metrics, however the optimal value of susceptance that maximises $\tau$ was different to the value that maximises both $\tau_\mathrm{H}$ and $\tau_\mathrm{A}$. It was found that true CCT would be improved by $47\%$ if the optimum susceptance loading for load C as measured by the analytic CCT was used instead of the original value of susceptance for this load. Given the quick assessment provided by an analytic CCT, an exhaustive study of all three-phase to ground faults on a network with their respective post-fault clearing strategies can be performed under a continuous range of loading conditions without having to do any formal fault study. Future research would then be required to test mathematical optimisation methods designed to find optimal loading distributions that improve the stability of the power system as measured by the CCT.

Despite its drawbacks, our analytic stability metric has the potential to inform optimal fault management strategies to improve system stability through parameteric investigation. Its key advantage is that it can be computed instantly once all the system parameter values for the pre-fault, fault-on and post-fault systems have been collected. This feature of stability metrics could be of use due to system dynamics becoming more unpredictable from the changing nature of loads \cite{yamashita2012modelling} and generation \cite{urdal2014system,Nguyen2008,Gautam2009} under the constraint of limited power flow through transmission lines. Particularly, investigations on the effect of low inertia on power system stability \cite{Ulbig,Tielens2012} can potentially benefit and this is an area of future work. Furthermore, optimisation techniques could be applied to analytic metrics to find regions in parameter space that increase power system stability in terms of CCT.


%





\ifCLASSOPTIONcaptionsoff
  \newpage
\fi



\bibliographystyle{IEEEtran}
\bibliography{ALLpapersAcademic,ALLbooksAndLectureNotes}
\end{document}